\documentclass[a4paper]{elsart}

\usepackage{epsfig} 
\usepackage{amsmath,amssymb,amsfonts} 
\usepackage{subfigure}
\usepackage{graphicx}          
                              
\journal{Automatica}
\usepackage{psfrag}

\newtheorem{definition}{Definition}

\newcommand{\R}{{\mathbb R}}
\newcommand{\abs}[1]{\left\vert #1 \right\vert}
\newcommand{\norm}[1]{\left\Vert #1 \right\Vert}

\newcommand{\sM}{\mathcal{M}}

\begin{document}

\begin{frontmatter}

\title{Contraction analysis of switched systems: the case of Caratheodory Systems and Networks} 

\author[]{Mario di Bernardo$^{a,b}$}
\author[First]{Davide Liuzza} 
\author[First]{Giovanni Russo}

\address[First]{Department of Systems and Computer Engineering,\\ University of Naples Federico II, 80125, Naples, Italy.\\ {\tt \{mario.dibernardo, davide.liuzza, giovanni.russo2\}@unina.it}}                                              
\address[Second]{Department of Engineering Mathematics, University of Bristol, BS8 1TR, Bristol, U.K.  {\tt \{m.dibernardo@bristol.ac.uk\}}}

\begin{abstract}                          
In this paper we extend to a generic class of piecewise smooth dynamical systems a fundamental tool for the analysis of convergence of smooth dynamical systems: contraction theory. We focus on switched systems satisfying Caratheodory conditions for the existence and unicity of a solution. After generalizing the classical definition of contraction to this class of dynamical systems,  we give sufficient conditions for global exponential convergence of their trajectories. The theoretical results are then applied to solve a set of representative problems including  proving global asymptotic stability of switched linear systems, giving conditions for incremental stability of piecewise smooth systems, and analyzing the convergence of networked switched linear systems.
\end{abstract}

\end{frontmatter}

\section{Introduction}
Piecewise-smooth dynamical systems are commonly used in Nonlinear Dynamics and Control to model devices of interest and/or synthesize discontinuous control actions e.g., \cite{Co:08}, \cite{BeBu:08}. Despite the large number of available results on their well-posedness and stability, there are few papers in the literature where the problem of assessing their incremental stability and convergence properties is discussed. 

The available results usually refer to specific classes of systems and rely on ad hoc continuity assumptions. For example, the problem of proving convergence for piecewise affine  continuous systems and networks is addressed in \cite{PaPo:05}, \cite{Pav_Wou_Nij_06}, \cite{PaPo:07}, \cite{PaWo:08} using a Lyapunov based approach. The methodology extends the approach of Demidovich for smooth dynamical systems expounded in \cite{Pav_Pog_Wou_Nij}. 

Contraction theory is an alternative approach used to study convergence between trajectories in smooth dynamical systems (see \cite{Loh_Slo_98} and references therein). The contraction approach is based on finding some metric under which the matrix measure of the system Jacobian can be proved to be definite negative over some convex set of phase space of interest. It can be shown that Demidovich approach is related to proving contraction using Euclidean norms and matrix measures. Nevertheless, using contraction it suffices to find some measure to study the Jacobian properties including non-Euclidean ones (e.g., $\mu_1$, $\mu_\infty$ etc).

Contraction theory has been used in a wide range of applications. For example, it has been shown that contraction is an extremely useful property to analyze coordination problems in networked control systems such as the emergence of synchronization or consensus 
(\cite{Loh_Slo_98,Pha_Slo_07,Rus_diB_Slo_09,Wan_Slo_05,Rus_diB_09b,Rus_diB_d,Rus_Slo_10}). Indeed, all trajectories of a contracting system can be shown to exponentially converge towards each other asymptotically. Therefore as shown in \cite{Wan_Slo_05}, this property can be effectively exploited to give conditions for the synchronization of a network of dynamical systems of interest. Recently, it has also been shown that the use of non Euclidean matrix measures can be used to construct an algorithmic approach to prove contraction \cite{Rus_diB_Slo_09} and to prove efficiently convergence in biological networks \cite{Rus_diB_Son_09}.

Historically,ideas closely related to contraction can be traced back
to~\cite{Hartmann} and even to~\cite{Lewis} (see
also~\cite{Pav_Pog_Wou_Nij}, ~\cite{Ang_02}, and e.g.~\cite{pde}, \cite{Jou_05} for a
more exhaustive list of related references).  For autonomous systems
and with constant metrics, the basic nonlinear contraction result
reduces to Krasovskii's theorem (\cite{SloLi}) in the
continous-time case, and to the contraction mapping theorem in the
discrete-time case (\cite{Loh_Slo_98}, \cite{Ber_Tsi_89}).

Despite the usefulness of contraction theory in applications, there is no consistent extension of this approach to the large class of piecewise smooth and switched dynamical systems. In \cite{Loh_Slo_00}, it is conjectured that, for certain classes of piecewise-smooth systems, contraction of each individual mode is sufficient to guarantee convergence of all the system trajectories towards each other, i.e. contraction of the overall system of interest. Also, in \cite{Rif_Slo_06}, it is noted that contraction theory can be extended to a class of hybrid systems under certain assumptions on the properties of the reset maps and switching signals.

The aim of this paper is to start addressing systematically the extension of contraction theory to generic classes of switched systems. Here, we focus on two types of switched systems of relevance in applications: (i) piecewise-smooth continuous (PWSC) systems (a class of state-dependent switched systems), and (ii) time-dependent switched systems (TSS). The goal is to obtain a set of sufficient conditions guaranteeing global exponential convergence of their trajectories. The extension of contraction theory to Filippov vector fields will be presented in a separate paper currently in preparation.

From a methodological viewpoint, to investigate contraction properties of PWSC and TSS we focus on systems whose vector fields satisfy Caratheodory conditions for the existence and uniqueness of an absolutely continuous solution (see \cite{Fi:88} and \cite{Co:08} for further details). We prove that, as conjectured in \cite{Loh_Slo_00}, for this class of systems, contraction of each individual mode suffices to guarantee convergence of all the system trajectories towards each other, i.e. contraction of the overall system of interest. Contrary to the results presented in  \cite{PaPo:05}, \cite{Pav_Wou_Nij_06}, \cite{PaPo:07}, \cite{PaWo:08}, we do not require finding incremental Lyapunov functions for the system of interest. We use instead, as commonly done for analysing contraction in smooth systems, a generic condition on the existence of some metric in which the Jacobian of each mode of the system of interest is definite negative.
We then apply the theoretical results to study convergence of some representative examples, including a network of time-switched systems.

A preliminary version of some of the results presented in this paper were presented at the $18^{th}$ IFAC World Congress in Milan in September 2011~\cite{RuDi:11}.

\section{Mathematical Preliminaries}
\label{sec:mathematicalpreliminaries}

Before presenting the main results of the paper, we introduce here some essential definitions and notation that will be used in the rest of the paper.

Let $x$ be an $n$-dimensional vector. We denote with $|x|$ the norm of the vector.  Let $A$ be a (real) matrix. Then, $\norm{A}$ denotes the norm of $A$.
We recall (see for instance~\cite{Mi_Li_07}) that, given a vector norm on Euclidean space ($\abs{\cdot}$), with its induced matrix norm $\norm{A}$, the associated \emph{matrix measure} $\mu$ is defined as the directional derivative of the matrix norm, that is,
\[
\mu(A) \,:=\;
\lim_{h \searrow 0} \frac{1}{h} \left(\norm{I+hA}-1\right).
\]
For example, if $\abs{\cdot}$ is the standard Euclidean 2-norm, then
$\mu(A)$ is the maximum eigenvalue of the symmetric part of $A$.
As we shall see, however, different norms will be useful for our applications.
Matrix measures are also known as ``\emph{logarithmic norms}'', a concept
independently introduced by Germund Dahlquist and Sergei Lozinskii in 1959,
\cite{dahlquist,lozinskii}.
The limit is known to exist, and the convergence is monotonic,
see ~\cite{strom,dahlquist}.

In what follows we report the analytic expression of some matrix measures used in the paper:
\begin{itemize}
\item $\mu_1(A)=\max_{j}\left(a_{jj}+\sum_{i\neq j}|a_{ij}|\right)$;
\item $\mu_2(A)=\max_{i}\lambda_i\left(\frac{1}{2}\left(A+A^T\right)\right)$;
\item $\mu_{\infty}(A)=\max_{i}\left(a_{ii}+\sum_{j\neq i}|a_{ij}|\right)$
\end{itemize}
More generally,  we will also make use of matrix measures induced  by  weighted vector norms, say $\abs{x}_{\Theta,i}= \abs{\Theta x}_i$, with $\Theta$ a constant
invertible matrix and $i=\{1,2,\infty\}$. Such measures, denoted with
$\mu_{\Theta,i}$, can be computed by using the following property: $\mu_{\Theta,i}(A)= \mu_i\left(\Theta A \Theta^{-1}\right)$, $\forall i=\{1,2,\infty\}$. 
Obviously, any other measure can be used.

%

We will also use the following definitions.
\begin{definition}[$K$-reachable sets \cite{Rus_diB_Son_09}]\label{defn:K-reachable}
Let $K>0$ be any positive real number.
A subset $C\subset\R^n$ is \emph{$K$-reachable} if, for any two points $x_0$ and
$y_0$ in $C$ there is some continuously differentiable curve 
$\gamma : \left[  0, 1 \right] \rightarrow C$ 
such that:
\begin{enumerate}
\item
$\gamma \left(0\right) = x_0$, 
\item
$\gamma\left(1\right)=y_0$ and 
\item
$\abs{\gamma '  \left(r\right)} \le K \abs{y_0 - x_0}$, $\forall r$. 
\end{enumerate}
\end{definition}
For convex sets $C$, we may pick $\gamma(r)=x_0+r(y_0-x_0)$, so 
$\gamma'(r)= y_0-x_0$ and we can take $K=1$.  Thus, convex sets are
$1$-reachable, and it is easy to show that the converse holds as well.

\begin{definition}[Flow of a Dynamical System]
Given a dynamical system 
%

\[
\dot{x}=f(x,t),\qquad x\in\mathcal{D}\subseteq\R^n, t\in\R,
\]

\noindent we define its flow $\varphi(s,t_0,\chi): \R^+ \times \R^+ \times \R^n \rightarrow \R^n$ as the operator (see \cite{mct} for a definition) such that 

\[
\frac{\partial}{\partial s}\varphi(s,t_0,\chi)=f(\varphi(s,t_0,\chi),t), \qquad \varphi(0,t_0,\chi)=\chi.
\] 
\end{definition}

In applications of the theory, it is often the case that $\mathcal{D}$ will
be a closed set, for example delimited by some hyperplane in the phase space, which could e.g. model constraints on the state variables of the system.
We remark here that for a non-open set, differentiability in $x$ means that the vector field $f(\bullet,t,)$ can be extended as a differentiable function to some open set which includes $\mathcal{D}$, and any continuity hypotheses with respect to 
$(t,x)$ hold on this open set.

\subsection{Caratheodory Solutions}\label{subsec:caratheodorysolutions}
We give the following preliminary definitions:

\begin{definition}
A function $g(t): [t_{0},+\infty] \rightarrow \R^{n}$ is said to be measurable if, for any real number $\alpha$, the set $\{t \in [t_{0},+\infty] :  g(t)>\alpha\}$ is measurable in the sense of Lebesgue.  
\end{definition}
\begin{definition}
A function $l(t): [t_{0},+\infty] \rightarrow \R$ is summable if the Lebesgue integral of the absolute value of $l(t)$ exists and is finite.
\end{definition}
\begin{definition}
A function $z(t) :[a,b] \rightarrow \R^{n}$ is absolutely continuous if for all $\varepsilon > 0$ there exists $\delta >0$ such that for each finite collection $[a_{1},b_{1}] \ldots [a_{n},b_{n}]$ of disjoint sets in $[a,b]$,  it holds that
$\sum_{k}\vert b_{k}-a_{k}\vert < \delta \implies \sum_{k} \vert z(b_{k})-z(a_{k}) \vert < \varepsilon$.
\end{definition}

Now we are able to give the definition  of a Caratheodory solution of a differential equation (see also \cite{Co:08} and references therein):

\begin{definition}[Caratheodory solutions]
Let us consider a domain $\mathcal{D}\subseteq\R^n$ and a dynamical system of the form:
\begin{equation}\label{eqn:PWS_NL}
\dot{x}(t) = f(x(t),t), \qquad x(t_0)=x_0,
\end{equation}

\noindent where $f: \mathcal{D} \times \R^{+} \rightarrow \R^{n}$. 
A Caratheodory  solution  for this system is an absolutely continuous function $x(t)$ that satisfies (\ref{eqn:PWS_NL}) for almost all $t\in [t_0, t_1]$ (in the sense of Lesbesgue), where $[t_0, t_1]$  is an interval where the solution $x(t)$ is defined. 
That is, a Caratheodory solution of (\ref{eqn:PWS_NL}) is an absolutely continuous function $x(t)$ such that:

\[
x(t) = x(t_0) + \int_{t_0}^t f(x(\tau),\tau)d\tau,\qquad t\in [t_0, t_1].
\]

In the common cases where the solution $x(t)$ can be extended forward in time, the same definition holds for each $t_1\geq t_0$.
\end{definition}

An useful result provides sufficient conditions for the existence of a Caratheodory solution of the system (\ref{eqn:PWS_NL}). 

\begin{thm}[Existence and Uniqueness of Caratheodory solutions]\label{thm:caratheodory}
A Caratheodory solution of system (\ref{eqn:PWS_NL}) exists if: 
\begin{enumerate}
\item for almost all $t \in [0,\infty]$, the function $x \mapsto f(x,t)$ is continuous for all $x \in \mathcal{D}$;
\item for each $x\in \mathcal{D}$, the function $t \mapsto f(x,t)$ is measurable in $t$;
\item for all $(x,t) \in \mathcal{D} \times [0, +\infty]$, there exist $\delta>0$ and  a summable function $m(t)$  defined on the interval $[t,t+ \delta]$ such that $\abs{f(x,t)}\le m(t)$.
\end{enumerate}

Moreover, the solution is unique, if the following additional condition is satisfied: 
\begin{itemize}
\item[(4)] $(x-y)^T (f(x,t)-f(y,t)) \le l(t) (x-y)^T(x-y)$, where $l(t)$ is a summable function.
\end{itemize}
\end{thm}

Notice that, as discussed in \cite{Fi:88}, p. 10 and proved in \cite{anbo:58}, equations that satisfy the  above (Caratheodory) conditions and those required for the uniqueness of a solution show continuous dependence on initial conditions.

\subsection{Caratheodory Systems}\label{subsec:pwscsystems}

Following \cite{BeBu:08} p.73,   we define a piecewise smooth dynamical system as follows.

\begin{definition}\label{def:pws}
Given a finite collection of disjoint, open and non empty sets $\mathcal{S}_1,\dots \mathcal{S}_p$ such that $\mathcal{D}=\bigcup_{i=1}^p\bar{\mathcal{S}}_i\subseteq\R^{n}$ is a connected set, a dynamical system $\dot{x}=f(x,t)$ is called a  {\em piecewise smooth dynamical system} (PWS) if it is defined by a finite set of ODEs

\[
f(x,t)=F_i(x,t)\qquad x\in\mathcal{S}_i.
\]

The intersection $\Sigma_{ij}:=\bar{\mathcal{S}}_i\cap\bar{\mathcal{S}}_j$ is either a lower dimensional manifold or it is the empty set. Each vector field $F_i(x,t)$ is smooth in both the state $x$ and the time $t$ for any $x\in\mathcal{S}_i$. Furthermore it is continuously extended on the boundary $\partial \mathcal{S}_i$. 
\end{definition}

We now introduce two important classes of switched systems that will be analyzed in the paper and for which a Caratheodory solution exists.
\newline
A piecewise smooth dynamical system is said to be continuous (PWSC) if the following two conditions hold:

\begin{enumerate}
\item the function $(x,t)\mapsto f(x,t)$ is continuous for all $x\in\R^n$ and  for all $t\geq t_0$;
\item the function $x\mapsto F_i(x,t)$ is continuously differentiable for  all $x\in\mathcal{S}_i$ and for all $t\geq t_0$. Furthermore the Jacobians $\frac{\partial F_i}{\partial x}(x,t)$ can be continuously extended on the boundary $\partial \mathcal{S}_i$.
\end{enumerate}

Notice that in order for condition (1) to be satisfied the functions $(x,t)\mapsto F_i(x,t)$ must be continuous for all $t\geq t_0$ and $x\in\mathcal{S}_i$, and, for all $x\in\Sigma_{ij}\neq\emptyset$ and all $t\geq t_0$, it must hold  $F_i(x,t)=F_j(x,t)$.


We now give the definition of a time-dependent switching system according to \cite{Li:03}.

\begin{definition}
A {\em time-dependent switching system} is a dynamical system of the form
\begin{equation}\label{eq:timeswitching}
\dot x = f(x,t,\sigma), \qquad x\in\R^n,  
\end{equation}
\noindent where $\sigma(t):[0,+\infty)\rightarrow \Sigma=\{1,2,\dots,p\}$ is a piecewise continuous time-dependent switching signal taking one over $p$ finite possible values.
\end{definition}

Note that according to this definition, we are excluding the case where infinite switchings occur over finite time so that  Zeno behavior cannot occur (see \cite{Li:03} for further details).   

We will refer to {\em Caratheodory systems} to indicate any PWS or TSS satisfying the conditions for the existence of a Caratheodory solution given in Theorem~\ref{thm:caratheodory}.
 
\section{Contraction Theory: a brief overview}
In this section, the notion of contraction is briefly summarized for a generic nonlinear system of the form:
\begin{equation}
\label{equ:main}
\dot x = h(x,t)
\end{equation}
where $h:\mathcal{D} \times \R \mapsto \R^n$ is an $m$-dimensional vector field assumed to be sufficiently {\em smooth}. (For further details see \cite{Loh_Slo_98,Rus_diB_Slo_09,Rus_diB_Son_09}.)

\begin{definition}
The smooth dynamical system (\ref{equ:main}) is said to be
{\em (infinitesimal) contracting} over a $K$-reachable set $\mathcal{C} \subseteq \R^n$, if there exists some matrix measure, $\mu$, and a positive scalar $c$, such that 
$$\mu\left[\frac{\partial h}{\partial x}(x,t)\right]  \le  - c,\qquad  \forall x\in \mathcal{C},  \forall t \ge 0.$$ 
The positive scalar $c$ is said to be the {\em contraction rate} of the system.
\end{definition}

The basic result of nonlinear contraction analysis states that, if a system is contracting, then all of its trajectories exponentially converge towards each other, see \cite{Loh_Slo_98,Rus_diB_Son_09}.

\begin{thm}[Contraction]
\label{theorem:contraction}
 Assume that (\ref{equ:main}) is contracting and let $\bar x(t)$ and $\tilde x(t)$ be any two of its solutions with initial conditions $\bar x(t_0)= \bar x_0\in \mathcal{C}$ and $\tilde x(t_0) = \tilde x_0\in\mathcal{C}$. Then, for any $t\ge t_0$, it holds that
 $$
 \abs{\bar x(t) - \tilde x(t)} \le K \abs{\bar x_0 - \tilde x_0} e ^{-ct\left(t-t_0\right)} 
 $$
\end{thm}

\section{Contraction of Caratheodory Systems}
Contraction theory has been properly studied mostly in the case of smooth nonlinear vector fields. The case of switched and hybrid systems is only marginally addressed in the existing literature \cite{Loh_Slo_00}, \cite{Rif_Slo_06}.
\newline
In this section, we seek sufficient conditions for the convergence of trajectories of PWSC systems (a generic class of systems with state-dependent switchings) and time-dependent switched systems. For the sake of clarity, we keep the derivation for the two cases separate.

\subsection{Contraction of PWSC systems}
We start with PWSC systems as defined in \S \ref{subsec:pwscsystems}. We can state the following result:

\begin{thm}\label{thm:contraction_pwsc}
 
Let $\mathcal{C}\subseteq\mathcal{D}$ be a $K$-reachable set.  Consider a generic PWS system of the form

\begin{equation}\label{eq:pwsc_esteso}
\dot{x}=f(x,t)=\begin{cases}
F_1(x,t) & x\in \mathcal{S}_1\\
\vdots & \\
F_p(x,t) & x\in \mathcal{S}_p
\end{cases}
\end{equation}

\noindent defined as in Definition \ref{def:pws} for all $x\in\mathcal{C}$ and with $\Sigma_{ij}$ smooth manifolds for all $i,j=1,\dots,p$. Suppose that:
\begin{enumerate}
\item it fulfills conditions for the existence and uniqueness of a Caratheodory solution given in \S \ref{subsec:caratheodorysolutions};
\item there exist a unique matrix measure such that
\[
\mu\left(\frac{\partial F_i}{\partial x}(x,t)\right) \le -c_i,
\]
\noindent for  all $x \in \bar{\mathcal{S}}_i$ and all $t\geq t_0$,  with $c_i$ belonging to a  set of positive scalars (in what follows, we will define $c := \min_{i}{c_{i}}$).
\end{enumerate}

\noindent Then, for every two solutions $x(t)=\varphi(t-t_0,t_0,\xi)$ and $y(t)=\varphi(t-t_0,t_0,\zeta)$ with $\xi,\zeta\in \mathcal{C}$, it holds that:
\[
|x(t)-y(t)|\leq K e^{-c(t-t_0)}|\xi-\zeta|,
\]
for all $t \geq t_0$ such that $x(t),y(t) \in \mathcal{C}$. If $\mathcal{C}$ is forward-invariant then all trajectories rooted in $\mathcal{C}$ converge exponentially towards each other. 

%

\end{thm}

{\bf Proof.}\ 
Given two points $x(t_0)=\xi$ and $y(t_0)=\zeta$ and a smooth curve $\gamma:[0,1]\mapsto \mathcal{C}$ such that $\gamma(0)=\xi$ and $\gamma(1)=\zeta$, we can consider $\psi(t,r):=\varphi(t-t_0,t_0,\gamma(r))$ as the solution of (\ref{eq:pwsc_esteso}) rooted in $\psi(t_0,r) = \gamma(r)$, with $r \in [0,1]$. Notice that $\psi(t,r)$ is continuous with respect to $r$ for  all $t$. Notice also that $\gamma$ can be chosen so that $\psi(t,r)$ is differentiable with respect to $r$ for almost all the pairs $(t,r)$. Let 

\begin{equation}\label{eq:w}
w(t,r) := \frac{\partial \psi}{\partial r}, \qquad \textit{a.e. in}\ t, \textit{ a.e. in}\ r.
\end{equation}

Thus we have:
$$
\frac{\partial w}{\partial t} = \frac{\partial}{\partial t}\left(\frac{\partial \psi}{\partial r}\right) = \frac{\partial}{\partial r}\left(\frac{\partial \psi}{\partial t}\right) = \frac{\partial }{\partial r}\left(f(\psi(t,r),t)\right), \qquad \textit{  a.e. in}\ t, \textit{  a.e. in}\ r,
$$

In what follows we will  use the shorthand notation \textit{a.e.} to denote the validity of a given expression almost everywhere in both $t$ and $r$, unless stated otherwise.  

Since 
$$
\frac{\partial }{\partial r}\left(f(\psi(t,r),t)\right) = \frac{\partial }{\partial x}f(\psi(x,t),t)\frac{\partial \psi(t,r)}{\partial r},\qquad a.e.,
$$
\noindent we can write:
\begin{equation}\label{eq:dinamica_w}
\frac{\partial}{\partial t} w(t,r) = A(\psi(t,r),t)w(t,r), \qquad a.e.,
\end{equation}

\noindent where we have denoted by $A(x,t)$ the Jacobian of the PWSC system (\ref{eq:pwsc_esteso}), which can be defined as:
\[
A(x,t)=\frac{\partial f}{\partial x}(x,t)=\begin{cases}
\frac{\partial F_1}{\partial x}(x,t), &  \forall x\in\mathcal{S}_1,\\
\vdots & \\
\frac{\partial F_p}{\partial x}(x,t), &  \forall x\in\mathcal{S}_p,
\end{cases}
\] 
for almost all the pairs $(x,t)$ apart from those points where $x \in \Sigma_{ij}$, for some $i,j$.


The next step is to show that the solution $t\mapsto w(t,r)$ of (\ref{eq:dinamica_w}) is a continuous function for any fixed $r \in [0,1]$.


Indeed, without loss of generality, consider the image of the curve $\gamma$ under the action of the flow $\varphi$ for a time $T$ such that the system trajectory rooted in $\gamma$ has either crossed the boundary once or it has not (in the case there are multiple switchings between $t_0$ and $T$, the same reasoning can be iterated). Furthermore, let us call $\tau(r)\in ]t_0, T[$ the time instant at which the trajectory eventually crosses the boundary. Suppose, without loss of generality, that at $t=\tau(r)$, the flow switches from region $\mathcal{S}_1$ to region $\mathcal{S}_2$. 
Then, we have:

\[
\psi(t,r)=\begin{cases}
\varphi_1(t-t_0,t_0,\psi(t_0,r)) & t_0\leq t < \tau(r), \\
\varphi_2(t-\tau(r),\tau(r),\varphi_1(\tau(r)-t_0,t_0,\psi(t_0,r))) & \tau(r)< t \leq T.  
\end{cases}
\]

%
Now, to show continuity of $w(t,r)$ with respect to time, from \eqref{eq:w} we need to evaluate the derivative of $\psi(t,r)$ over the interval $]t_0,T[$. We have:
\begin{equation}\label{eq:derivata_psi}
\frac{\partial \psi}{\partial r}(t,r)=\begin{cases}
\frac{\partial}{\partial r}[\varphi_1(t-t_0,t_0,\psi(t_0,r))] & t_0\leq t < \tau(r), \\
\frac{\partial}{\partial r}[\varphi_2(t-\tau(r),\tau(r),\varphi_1(\tau(r)-t_0,t_0,\psi(t_0,r)))] & \tau(r)< t \leq T,  
\end{cases}
\end{equation}
Continuity of $w(t,r)$ is then guaranteed if
\begin{equation}\label{eq:limite_continuita_w}
\lim_{t\rightarrow \tau(r)^{-}}\frac{\partial}{\partial r}[\varphi_1(t-t_0,t_0,\psi(t_0,r))]=
\lim_{t\rightarrow \tau(r)^{+}}\frac{\partial}{\partial r}[\varphi_2(t-t_0-\tau(r),\tau(r),\varphi_1(\tau(r)-t_0,t_0,\psi(t_0,r)))],
\end{equation}

We have
\begin{equation}
\frac{\partial}{\partial r} \varphi_1(s,t_0,\chi)=\frac{\partial \varphi_1}{\partial \chi} \frac{\partial \chi}{\partial r}
\end{equation}
with $s:=t-t_0$ and $\chi=\psi(t_0,r) := \psi^0$. Hence, taking the limit $t \rightarrow \tau(r)^-$, the left-hand side of (\ref{eq:limite_continuita_w}) can be written as:
\begin{equation} \label{eq:lhs}
\frac{\partial \varphi_1}{\partial \psi^0}(\tau(r)-t_0,t_0, \psi^0)\frac{\partial \psi^0}{\partial r}
\end{equation}

Also 
$$
\frac{\partial}{\partial r} \varphi_2(s(t,r),\hat t_0(r),\chi(r)) = \frac{\partial \varphi_2}{\partial s}\frac{\partial s}{\partial r} + \frac{\partial \varphi_2}{\partial \hat t_0} \frac{\partial \hat t_0}{\partial r}+\frac{\partial \varphi_2}{\partial \chi}\frac{\partial \chi}{\partial r}
$$
where
\begin{eqnarray}
s(t,r) &:=& t-\tau(r) \\
\hat t_0(r) &:=& \tau(r) \\
\chi(r) &:=& \varphi_1\left(-s(t_0,r),t_0,\psi^0\right)
\end{eqnarray}

Now, we observe that 
\begin{eqnarray}
\frac{\partial \varphi_2}{\partial s} &=& F_2(\varphi_2(s(t,r),\tau(r),\chi(r)),t) \label{eq:f2s}\\
\frac{\partial \chi}{\partial r} &=& \frac{\partial \varphi_1}{\partial s}\tau^\prime(r)+\frac{\partial \varphi_1}{\partial \psi^0}\frac{\partial \psi^0}{\partial r}
\end{eqnarray}
where
$$
\frac{\partial \varphi_1}{\partial s}=F_1(\varphi_1(-s(t_0,r),t_0,\psi^0),t)
$$

Taking the limit for $t \rightarrow \tau(r)^\pm$, we have:
$$
\lim_{t \rightarrow \tau(r)^+} s(t,r)=0
$$
hence, since
\begin{equation} \label{eq:id}
\varphi_2(0, \tau(r), \chi(r), t)= \chi(r)
\end{equation}
we then obtain that (\ref{eq:f2s}) yields in the limit
$$
\frac{\partial \varphi_2}{\partial s}=F_2(\chi(r),\tau(r))=F_2(\varphi_1(-s(t_0,r),t_0,\psi^0),\tau(r))
$$

Moreover, in the same limit, $t \rightarrow \tau(r)$, from (\ref{eq:id}) we have:
$$
\frac{\partial \varphi_2}{\partial \hat t_0} = \frac{\partial \chi}{\partial \hat t_0}=0,
$$

\noindent and

$$
\frac{\partial \varphi_2}{\partial \chi} = \frac{\partial \chi}{\partial \chi}=I.
$$

Therefore, the right-hand side of (\ref{eq:limite_continuita_w}) can be written as
\begin{eqnarray} \label{eq:rhs1}
&-&F_2(\varphi_1(-s(t_0,r),t_0,\psi^0),\tau(r)) \tau^\prime(r) \nonumber \\ &+& F_1(\varphi_1(-s(t_0,r),t_0,\psi^0),\tau(r)) \tau^\prime(r)+ \frac{\partial \varphi_1}{\partial \psi^0}\frac{\partial \psi^0}{\partial r}
\end{eqnarray}

From the assumption that the system vector field is continuous when $t=\tau(r)$, continuity of $w(t,r)$ with respect to $t$ is then immediately established by comparing (\ref{eq:lhs}) and (\ref{eq:rhs1}).

Now, we turn again our attention to equation (\ref{eq:dinamica_w}).
Fixing $r$ to any value between 0 and 1, the Jacobian can be calculated and (\ref{eq:dinamica_w}) can be solved to obtain (in the sense of Lebesgue):
$$
\begin{array}{*{20}l}
w(t+h, r) = w(t, r) + \int_{t}^{t+h}A(\psi(\vartheta, r),\vartheta)w(\vartheta, r)d\vartheta =\\
= w(t, r) + A(\psi(t,r),t)w(t, r)h +\\ 
\int_{t}^{t+h}\left(A(\psi(\vartheta, r),\vartheta)w(\vartheta,r)-A(\psi(t,r),t)w(t,r)\right)d\vartheta, \qquad \textit{  a.e. } t,
\end{array} 
$$

\noindent with $h$ being a positive scalar.

Thus, from the above expression we have

\begin{eqnarray}
|w(t+h,r)|\leq & & \left|\left|I + hA(\psi(t,r),t)\right|\right|  \left|w(t, r)\right| \nonumber\\ &+& 
\int_{t}^{t+h}\left|A(\psi(\vartheta, r),\vartheta)w(\vartheta,r)-A(\psi(t,r),t)w(t,r)\right|d\vartheta 
\end{eqnarray}
\noindent Then, subtracting $|w(t, r)|$ from both sides of the equation and dividing by $h$ we obtain

$$
\begin{array}{*{20}l}
\frac{1}{h}\left(\abs{w(t+h, r)}-\abs{w(t,r)}\right)  \le \frac{1}{h}\left(\norm{I+hA(\psi(t,r),t)}-1\right)\abs{w(t,r)}+\\ \frac{1}{h}\int_{t}^{t+h} \abs{A(\psi(\vartheta, r),\vartheta)w(\vartheta,r)-A(\psi(t,r),t)w(t,r)}d\vartheta, \qquad \textit{  a.e. } t, 
\end{array} 
$$

Thus, taking the limit as $h\searrow 0$ yields:
$$
\frac{d}{dt}\abs{w(t, r)} \le -c \abs{w(t,r)}, \qquad a.e.
$$

Notice that the above expression holds for all those pairs $t$ and $r$ where the Jacobian $A(\cdot)$ is defined.
Let now $M(t):= -c(t-t_0)$, from the above expression it follows that:
$$
\frac{d}{dt}\left(\abs{w(t, r)}e^{-M(t)}\right)\le 0, \qquad a.e.
$$

Now, since $e^{-M(t)}$ is an increasing function and since the function $t\mapsto w(t,r)$ is continuous, the above inequality implies that:
$$
\abs{w(t, r)} \le \abs{w(t_0, r)}e^{-c(t-t_0)} \le K \abs{\xi - \zeta}e^{-c(t-t_0)}.
$$

As the function $\psi(t,r)$ is continuous and, for all $t$, the function $w(t,r)$ is defined for almost all $r$,  we have:
$$
\psi(t,1) - \psi(t,0) = \int_{0}^1w(t,s)ds.
$$
Thus:
$$
\abs{x(t)-y(t)} \le K \abs{\xi-\zeta}e^{-c(t-t_0)},
$$
and the theorem remains proved.

Obviously, if $\mathcal{C}$ is forward-invariant, then trajectories rooted in $\mathcal{C}$ will exponentially converge towards each other.

%
%
%
%
%
%
%
\subsection{Contraction of TSS}
The conditions used to prove contraction of PWSC systems can be immediately extended to generic systems affected by time-dependent switchings as detailed below.

\begin{thm}\label{thm:contraction_timeswitching}
Consider an invariant $K$-reachable set $\mathcal{C}\subseteq\R^n$ and a time-dependent switching system as in (\ref{eq:timeswitching}). Suppose that:
\begin{enumerate}
\item it fulfills conditions for the existence and uniqueness of a Caratheodory solution given in \S \ref{subsec:caratheodorysolutions};
\item the function $(x,t)\mapsto f(x,t,\sigma)$ is continuous for all $x\in\mathcal{C}$, for all $t\geq t_0$ and for all $\sigma\in\Sigma$;
\item the function $x\mapsto f(x,t,\sigma)$ is continuously differentiable for all $x\in\R^n$, for all $t\geq t_0$ and for all $\sigma\in\Sigma$;
\item there exist a unique matrix measure such that
\[
\mu\left(\frac{\partial f}{\partial x}(x,t,\sigma)\right) \le -c_\sigma^2,
\]
\noindent for all $x \in\mathcal{C}$, for  all $t\geq t_0$ and for all $\sigma\in\Sigma$, with $c_\sigma$ belonging to a  set of real scalars (in what follows, we will define $c^2 := \min_{\sigma \in \Sigma} c^2_\sigma$).
\end{enumerate}

Then, for every two solutions $x(t)=\varphi(t,t_0,\xi)$ and $y(t)=\varphi(t,t_0,\zeta)$ with $\xi,\zeta\in \mathcal{C}$, it holds that:

\[
|x(t)-y(t)|\leq K e^{-c^2(t-t_0)}|\xi-\zeta|,
\]
\end{thm}

{\bf Proof.}\ 
The proof follows similar steps to that of Theorem \ref{thm:contraction_pwsc}. In particular, given  points $x(t_0)=\xi$ and $y(t_0)=\zeta$ and a smooth curve $\gamma:[0,1]\mapsto \mathcal{C}$ such that $\gamma(0)=\xi$ and $\gamma(1)=\zeta$, we can consider $\psi(t,r):=\varphi(t-t_0,t_0,\gamma(r))$ as the solution of (\ref{eq:timeswitching}) rooted in $\psi(t_0,r) = \gamma(r)$, with $r \in [0,1]$. Let 

$$
w(t,r) := \frac{\partial \psi}{\partial r}, \qquad \textit{  a.e. } t.
$$ 

As in the proof of Theorem \ref{thm:contraction_pwsc}, we can write:
$$
\frac{\partial w}{\partial t} = \frac{\partial}{\partial t}\left(\frac{\partial \psi}{\partial r}\right) = \frac{\partial}{\partial r}\left(\frac{\partial \psi}{\partial t}\right) = \frac{\partial }{\partial r}\left(f(\psi(t,r),t,\sigma)\right), \qquad \textit{  a.e. } t,
$$

\noindent and then

\[
\frac{\partial}{\partial t} w(t,r) = A(\psi(t,r),t,\sigma)w(t,r), \textit{  a.e. } t,
\]

\noindent with $A=\frac{\partial}{\partial x}f(x,t,\sigma)$ being the Jacobian of $f(x,t,\sigma)$ for almost all $t\geq t_0$. 
Notice that, differently from the case of PWSC systems, here the Jacobian is discontinuous only with respect to time $t$ due the fact that the function $\sigma(t)$ is piecewise constant.
However, we can show that the function $t\mapsto w(t,r)$ is continuous by considering again (\ref{eq:lhs}) and (\ref{eq:rhs1}). In this case, the switching instant $\tau(r)$ is independent from $r$ and therefore all terms containing $\tau^\prime(r)$ in (\ref{eq:rhs1}) cancel out. The equality of (\ref{eq:lhs}) and (\ref{eq:rhs1}) then immediately follows and the rest of proof becomes identical to that of Theorem \ref{thm:contraction_pwsc}.

\subsection*{Remarks}
\begin{itemize}
\item Our results on the contraction of PWSC systems can be interpreted following the approach presented in \cite{CaFr:06} where the asymptotic stability of piecewise linear systems obtained by the continuous matching of two stable linear systems is discussed. Specifically, under the conditions of Theorem~\ref{thm:contraction_pwsc}, we can state that {\em the continuous matching of any number of nonlinear contracting vector fields is also contracting}. As in the case of asymptotic stability discussed in \cite{CaFr:06}, guaranteeing incremental stability of switched systems is not trivial, even when they are obtained by continuously matching contracting vector field. Thus, the sufficient conditions derived in this paper can be useful for the analysis of incremental stability in switched systems and the design of stabilizing switched control inputs. 

\item The results reported so far do not include the case of differential equations with discontinuous right-hand side or Filippov systems where so-called sliding motion is possible \cite{BeBu:08}. This is ongoing work that will be reported elsewhere \cite{future_paper}.
\end{itemize}

\section{Partial Contraction of PWSC and TSS}
In the previous section, contraction theory was extended to piecewise continuous and time-dependent systems. In both cases, we were able to show that the contraction still implies asymptotic convergence between trajectories.
Often in applications, it is desirable to prove (or certify) that, at steady state, all trajectories of a given system exhibit some property regardless of their initial conditions.

In the case of smooth dynamical systems, the concept of partial contraction was introduced in  \cite{Wan_Slo_05} to solve this problem. The idea is to introduce an appropriately constructed auxiliary or {\em virtual system}, embedding the solutions of the system of interest as its particular solutions. If the virtual system is proved to be contracting, then all of its solutions will converge towards a unique trajectory. In turn, this imples that all trajectories of the system of interest, embedded in the virtual system by construction, will also converge towards this solution.

The most notable application of partial contration is its use to prove convergence of trajectories of all nodes of a network of dynamical systems towards each other as for example is required in synchronization or coordination problems. In that case, the virtual system is constructed so that trajectories of the network nodes are its particular solutions. Proving contraction of the virtual system then implies convergence of all node trajectories towards the same synchronous evolution (see \cite{Wan_Slo_05}, \cite{Slo_Wan_Rif_98}, \cite{Rus_Slo_10}, \cite{Rus_diB_09b}, \cite{Rus_diB_Slo_09} for further details and applications).

Using the extension of contraction to switched Caratheodory systems presented in this paper, we can also extend partial contraction to this class of systems. For example, we can prove the following result for PWS systems.

\begin{thm}
\label{theorem:partial_contraction}
Consider system a PWS of the form (\ref{eq:pwsc_esteso}) and assume that there exists some system
\begin{equation}\label{eqn:virtual_system}
\dot y= v(y,x,t)
\end{equation}
such that:
\begin{itemize}
\item $v(x,x,t) = f(x,t)$;
\item $v(y,x,t)$ is contracting in the Caratheodory sense with respect to $y$ and for any $x$.
\end{itemize}
Let $\tilde y$ be the unique solution towards which all trajectories of (\ref{eqn:virtual_system}) converge. That is, there exists some $c>0$ such that, for any solution of (\ref{eqn:virtual_system}), say $y(t)$:
$$
\abs{y(t) - \tilde y(t)} \le \abs{y(t_0)-\tilde y(t_0)}e^{-c(t-t_0)}
$$
Then, all the solutions of (\ref{eq:pwsc_esteso}) converge towards $\tilde y$, i.e.
$$
\abs{x(t)-\tilde y(t)}\rightarrow 0, \ \ \ t \rightarrow + \infty
$$
\end{thm}
 System (\ref{eq:pwsc_esteso}) is said to be \emph{partially contracting} while system (\ref{eqn:virtual_system}) is termed as virtual system.
 
{\bf Proof.}
Indeed, we only need to observe that by construction any solution of (\ref{eq:pwsc_esteso}), say $x(t)$, is also a solution of the virtual system. Now, since (\ref{eqn:virtual_system}) is contracting, then all of its solutions will converge towards $\tilde y$. This in turn implies that
$$
\abs{x(t) - \tilde y(t)} \rightarrow 0 \ \ \ \ \ a.e.
$$
as $t\rightarrow + \infty$.
\qed

We remark here that Theorem \ref{theorem:partial_contraction} allows to prove that all the solutions of the real system of interest converge, at steady state, towards a unique solution even if it is not contracting. The key point of such a result is indeed that of constructing a contracting system which embeds the solutions of the real system. The main implication of this is that trajectories of the real system will converge towards each other but convergence will not be exponential and, in general, it may be non-uniform. That is, Theorem \ref{theorem:partial_contraction} only ensures that, after a sufficiently long time, distances between all the solutions of (\ref{eq:pwsc_esteso}) shrink. In some special case, see e.g. \cite{Rus_Slo_10}, the dimensionality of the virtual system is lower than that of the real system of interest: this is typically the case of systems with symmetries, such as \emph{Quorum Sensing} networks. A notable example of use of virtual system can be found in \cite{Rus_Slo_11}. We also remark that Theorem \ref{theorem:partial_contraction} can be straightforwardly extended to time-dependent switched systems. The proof follows exactly the same steps of those used to prove Theorem \ref{theorem:partial_contraction} and hence it is omitted here for the sake of brevity.

\subsubsection*{Example 1}
As an example illustrating the key features of partial contraction and virtual systems, consider a PWSC system of the form
\begin{equation}\label{eqn:example_virtual}
\dot x = L(t,x) x.
\end{equation}
Notice that such a system may e.g. model a networked control system or a network of biochemical reactions.

We assume that the system is not contracting. That is, the Jacobian matrix
$$
\frac{\partial L}{\partial x}x + L(t,x),
$$
does not have any uniformly negative matrix measure. We also assume that there exists a uniformly negative matrix measure for $L(t,x)$, i.e.
$$
\exists \ \mu \ : \ \mu\left(L(t,x)\right) \le -c, \qquad c > 0  \quad a.e.
$$
Clearly, in this case, system (\ref{eqn:example_virtual}) is not contracting nevertheless Theorem \ref{theorem:partial_contraction} can be used to show that, at steady state, all trajectories of (\ref{eqn:example_virtual}) converge towards a unique solution.  
In particular, consider the system
$$
\dot y = v(y,x,t) = L(t,x) y,
$$
where $x$, the state variable of the original system, is seen as an external input. It is straightforward to check  that
$$
v(x,x,t) = L(t,x) x = \dot x,
$$
and hence it is a virtual system in the sense of Theorem \ref{theorem:partial_contraction}. Moreover, the Jacobian matrix of the virtual system is simply
$$
J (t,x ) = \frac{\partial L(t,x) y}{\partial y} = L(t,x).
$$
Since we assumed that there exists a uniformly negative matrix measure for $L(t,x)$, the virtual system is contracting for any $x$. Therefore, all of its solutions will converge to a unique trajectory, say $y^\ast$, such that:
$$
\dot y^\ast = L(t,x) y^\ast
$$
Since the solutions of the real system are also particular solutions ofthe virtual system, it follows that 
$$
\abs{x(t) - y^\ast (t)}\rightarrow 0 \ \ \ a.e.
$$
That is, all solutions of the real system will also converge towards $y^\ast$ and, hence, towards each other.

\section{Applications}
The extension of contraction and partial contraction to Caratheodory systems is a flexible tool that can be used for analysis and design as illustrated by means of two representative examples described in this section.

\subsection{Stability of Piecewise Linear Systems}
Using the concept of contraction for PWS systems, it is possible to easily prove the following result to assess the stability of piecewise linear systems  of the form
\begin{equation} \label{eq:PWL}
\dot{x}=A(t,\sigma) x,
\end{equation}
where $x \in \mathbb{R}^n$ and $\sigma(t):[0,+\infty)\rightarrow \Sigma$ is the switching signal with $\Sigma$ being a finite index set.
Several stability results for this class of systems are available in the literature (see \cite{Lin_Ant_09} for an extensive survey). A classical approach is that of finding conditions on the (finite) set of matrices $A(t,\sigma)$ guaranteeing the existence of some common quadratic Lyapunov function (CQLF). In \cite{Lin_Ant_09} (Theorem 8, p. 311), it is proven that the origin is a globally asymptotically stable equilibrium of  (\ref{eq:PWL}) if and only if there exist a full column rank matrix  $L\in\R^{m\times n}$ with $m\geq n$ and a family of matrices $\bar{A}_i\in\R^{m\times m}$ such that $\mu_{\infty}\left(\bar{A}_i \right)<0$ for all $i=\sigma\in\Sigma$. 

Here we show that a related stability condition can be immediately obtained by applying contraction theory. Indeed, we can prove the following result.

\begin{cor}\label{cor:PWL}
Given a piecewise linear system of the form (\ref{eq:PWL}), if the matrices $A(t,\sigma)$ are bounded and measurable for any $\sigma$ and there exist some matrix measure such that
\begin{equation}\label{eqn:matrix_meas_neg}
\mu \left(A(t,\sigma)\right) \le - c, \quad c > 0 \quad \forall t\in \R^+,\quad\forall \sigma \in \Sigma,
\end{equation}
then, all solutions of (\ref{eq:PWL}) converge asymptotically towards the origin independently from the switching sequence. 
\end{cor}

{\bf Proof.} Under the hypotheses, system (\ref{eq:PWL}) satisfies Theorem \ref{theorem:partial_contraction} and therefore is contracting with all of its trajectories converging towards each other. Since, $x(t)=0$ is also a trajectory of (\ref{eq:PWL}), the proof immediately follows. 

As an example, take the system:
\begin{equation}
\dot{x}=A(\sigma)x, \quad \sigma \in \{1,2\}
\end{equation}
with 
\begin{equation*}
A(1)=\begin{pmatrix}
-1.0 & 1.5 \cr
0.8 & -3.0
\end{pmatrix}, \quad A(2)=\begin{pmatrix}
-3.0 & 1.0 \cr
2.0 & -1.5
\end{pmatrix}. 
\end{equation*}

Note that using the matrix measure $\mu_1$ induced by the $1$-norm, we have $\mu_1[A(1)]<0$ and $\mu_1[A(2)]<0$. Hence, it is immediate to prove asymptotic convergence of all solutions towards each other and onto the origin using Corollary \ref{cor:PWL}, for arbitrarily switching signal $\sigma(t)$.

We wish to emphasize that the result based on contraction embeds as a special case the stability condition reported in \cite{Lin_Ant_09}. Indeed, setting $L=I$ and $\bar A_i=A(\sigma)$ in Theorem 8, p. 311 in \cite{Lin_Ant_09} is equivalent to using Corollary \ref{cor:PWL} with the matrix measure $\mu_\infty$. Moreover, the proof based on contraction can also be extended to nonlinear switched systems.

\subsection{Incremental Stability of PWSC and TSS systems}
Results stated in this paper give us a powerful tool to easily show exponential incremental stability for switched Caratheodory systems. 
As an example, we consider a switched version of the biological system used in \cite{Rus_diB_Son_09} to illustrate the application of contraction to biochemical networks.
Specifically, following \cite{Rus_diB_Son_09}, we consider the externally driven transcriptional system described by the equations:
\begin{eqnarray}
\dot{x}&=& u(t)\left(X_{tot}-x-y\right)-\delta x +k_1y-k_2\left(e_T-y\right)x,\label{eq:transcriptional_smooth_xp}\\
\dot{y}&=& -k_1y+k_2\left(e_T-y\right)x. \label{eq:transcriptional_smooth_yp}
\end{eqnarray}

In the above model, $x$ represents the concentration of the transcription factor $X$, $y$  the concentration of the complex protein-promoter $Y$ with the constant $e_T$ being its total concentration. The external signal $u(t)$ represents the concentration of a transcriptional factor which inactivates transciptor $X_0$ (and consequently $X$) through a phosphorylation reaction, while the scalars $k_1$ and $k_2$ represent the binding and dissociation rates, respectively. Finally, $X_{tot}$ represents the total constant amount of reactant concentrations $X_{tot}=X+Y$.


We now assume that the smooth system \eqref{eq:transcriptional_smooth_xp}-\eqref{eq:transcriptional_smooth_yp} is also affected by the possible presence of an extra term $v(t)$, of the form:

\[
v(t)=\begin{cases}
0 & \textit{ if } x\leq h   \\
-\beta\left[x-h\right]  & \textit{ if } x> h   
\end{cases}.
\]

From the physical viewpoint, the term $v(t)$ models a degradation (see e.g. \cite{Sza_Ste_Per_06}) on the state variable $x$ which becomes active when the value is above a certain threshold, say $h$. That is, the degradation is activated when the total amount of $x$ is greater than $h$ and its effect can be modelled as an additive term.

Using as in \cite{Rus_diB_Son_09} the simple change of variables $x_t=x+y$, the resulting model can be rewritten as:

\begin{eqnarray}
\dot{x_t}&=& u(t)\left(X_{tot}-x_t\right)-\delta x_t +\delta y+\tilde v(t),\label{eq:transcriptional_nonsmooth_xp}\\
\dot{y}&=& -k_1y+k_2\left(e_T-y\right)\left(x_t-y\right), \label{eq:transcriptional_nonsmooth_yp}
\end{eqnarray}

\noindent where the discontinuous term $\tilde v(t)$ is given by:

\[
\tilde v(t)=\begin{cases}
0 & \textit{ if } x_t-y\leq h   \\
-\beta\left[x_t-y-h\right]  & \textit{ if } x_t-y> h   
\end{cases}.
\]

To prove contraction, and hence global incremental stability of this switched system, we need to derive the Jacobian which, in this case, is the discontinuous function:
\[
J(x,t)=\begin{cases}
J_s & \textit{ if }  x\leq h  \\
J_s+J_{ns} & \textit{ if }  x> h  
\end{cases},
\]

\noindent where 
\[
J_s=\left[
\begin{array}{cc}
-u(t)-\delta & \delta\\
k_2(e_T-y)\quad & -k_1+k_2(-e_T-x_t+2y)
\end{array}
\right],
\]

\noindent and

\[
J_{ns}=\left[
\begin{array}{rr}
-\beta  & \beta\\
 0 & \quad 0
\end{array}
\right].
\]

Using $\mu_{\infty}(\cdot)$ as a matrix measure, we find that $\mu_{\infty}(J_s)$ is negative if the following inequalities hold:

\begin{equation}\label{eq:inequality_jacobiano1}
-u(t)-\delta +|\delta|<-c_1;
\end{equation}

\begin{equation}\label{eq:inequality_jacobiano2}
-k_1+k_2(-e_T-x_t+2y)+\left|k_2(e_T-y)\right|<-c_2;
\end{equation}

\noindent for $c_1,c_2>0$.
\newline
As shown in \cite{Rus_diB_Son_09}, the first inequality is always satisfied as the system parameters and the periodic input are assumed to be positive. Furthermore, taking into account that, for physical reasons, the term $e_T-y\geq 0$, it is easy to prove that inequality \eqref{eq:inequality_jacobiano2} is also fulfilled.  

\begin{figure}
\includegraphics[width=14cm]{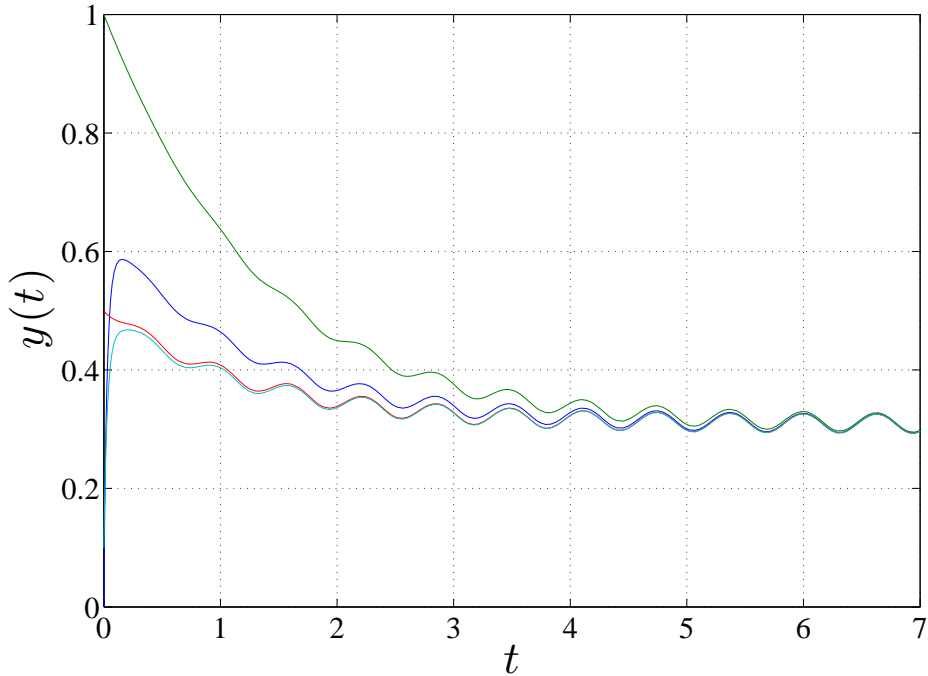}
\caption{Convergence of trajectories of the transcriptional module starting from different initial conditions towards the same unique periodic orbit. Simulation were carried out with the following parameter values: $k_1=0.5$, $k_2=5$, $X_{tot}=1$, $e_T=1$, $\delta=20$, $\beta=1$, $h=0.01$. The periodic input was set to $1.5+2\sin\left(10t\right)$.}%
\label{fig:sincronizzazione_y}
\end{figure}

We now have to consider the effect of the switching by looking at the measure of the matrix $J_s+J_{ns}$. It is immediate to see that, $\mu_\infty(J_s+J_{ns})$ is also negative if inequalities (\ref{eq:inequality_jacobiano1}) and (\ref{eq:inequality_jacobiano2}) are satisfied. 
Hence, according to Theorem~\ref{thm:contraction_pwsc}, the switched biochemical network under investigation is contracting and is therefore incrementally stable. This also implies that, as discussed in \cite{Rus_diB_Son_09}, the transcriptional network continues to be entrainable even when an additional discontinuous term is added to the model. This is confirmed by the numerical simulation reported in Fig.\ref{fig:sincronizzazione_y}.

\begin{figure}
\includegraphics[width=14cm]{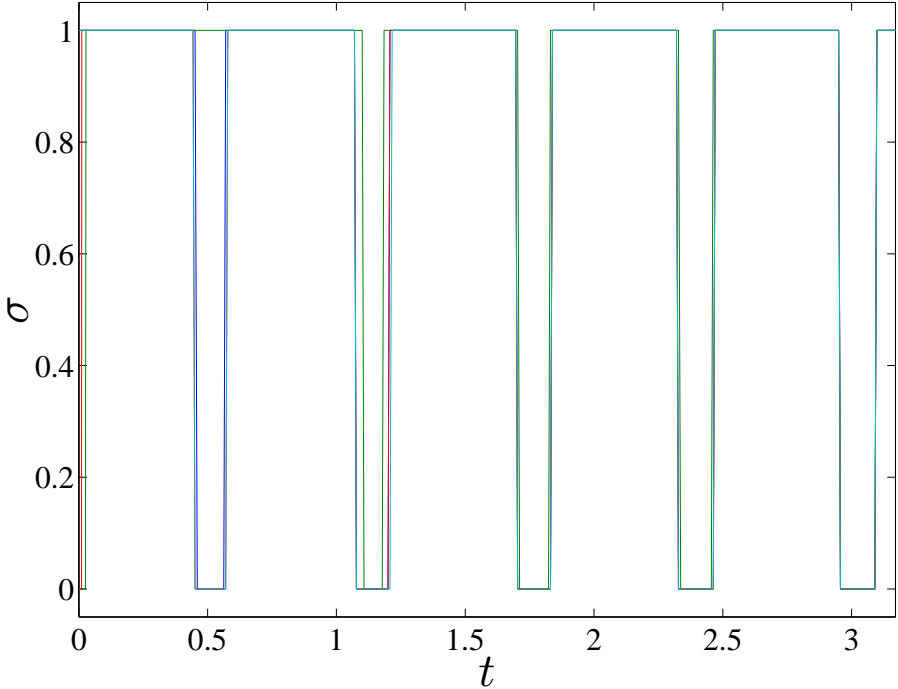}
\caption{Synchronization of the switching signals associated to trajectories of the modified transcriptional module starting from different initial conditions}%
\label{fig:sincronizzazione_sigma}
\end{figure}

Fig.\ref{fig:sincronizzazione_sigma} shows that, as expected, the switching signals associated to trajectories starting from different initial conditions also synchronize asymptotically. 

\subsection{Convergence of networks of time-switching systems}
Contraction analysis can be an invaluable tool to study convergence of networked systems as  proposed in \cite{Wan_Slo_05}. Here, we use the extension of contraction to Caratheodory systems presented in this paper to derive conditions guaranteeing convergence of a network of diffusively coupled switched linear systems.  Specifically, we consider a network of the form:
\begin{equation}\label{eqn:network_dyn}
\dot x_i = A(\sigma (t) )x_i + \Gamma \sum_{j \in N_i} \left[x_j - x_i \right],
\end{equation}
where $x_i \in \R^{n}$ represents the state vector of node $i$, $N_i$ denotes the set of the neighbors of the $i$-th  node whose cardinality  (i.e. the degree of the $i$-th network node) is denoted with $d_i$. In the above equation, $\Gamma$ is a coupling matrix, often termed as inner-coupling matrix in the literature. In what follows the eigenvalues of the network Laplacian matrix ($L$) are denoted with $\lambda_i$; $\lambda_2$ being the smallest non-zero Laplacian eigenvalue (algebraic connectivity). We assume $A(\sigma(t))$ to be bounded and measurable.

We will now show that, by using contraction,  a sufficient condition can be derived ensuring all the solutions of the network nodes globally exponentially converge, almost everywhere, towards the $n$-dimensional linear subspace\footnote{It is straightforward to check that this subspace is flow invariant for the network dynamics.} $\sM_s :=\left\{x_1=\ldots=x_N\right\}$. In what follows, we will denote by $s(t)$ the common asymptotic behavior of all nodes on $\sM_s$. Note that $s(t)$ is obviously a solution of each isolated node of (\ref{eqn:network_dyn}), i.e. $\dot s(t) = A(\sigma(t))s(t)$. We will also say that the network nodes are coordinated (or that the network is coordinated) if 
$$
\lim_{t \rightarrow \infty}{\abs{x_i(t) -s(t)}}= 0, \quad a.e.
$$
In the special case where $s(t)$ exhibits an oscillatory behavior, we will say that all network nodes are synchronized (or that the network is synchronized).

\begin{thm}\label{thm:linear_diffusive}
The trajectories of all nodes in the network (\ref{eqn:network_dyn}) exponentially converge  towards each other almost everywhere (i.e., the network is coordinated a.e.) if (i) the topology of the network is connected and (ii) there exist some matrix measure, $\mu$, such that:
$$
\mu\left(A(\sigma(t))- \lambda_2\Gamma\right) \leq -c, \quad c > 0,
$$
for all $\sigma \in \Sigma$ and for almost all $t$.
\end{thm}

Before presenting the proof of the Theorem, we report here two useful results, \cite{Arc_unp}.

\begin{lem}\label{lem:kronecker}
Let $\otimes$ denote the Kronecker product. The following properties hold:
\begin{itemize}
\item $\left(A \otimes B \right) \left(C \otimes D\right) = \left(A C\right) \otimes\left(BD\right)$; 
\item if $A$ and $B$ are invertible, then $\left(A \otimes B \right)^{-1} =A^{-1}\otimes B^{-1}$;
\end{itemize}
\end{lem}
\begin{lem}\label{lem:schur}
For any $n\times n$ real symmetric matrix, $A$, there exist an orthogonal $n\times n$ matrix, $Q$, such that
\begin{equation}\label{eqn:schur}
Q^TAQ= U,
\end{equation}
where $U$ is an $n\times n$ diagonal matrix.
\end{lem}

{\bf Proof (Theorem \ref{thm:linear_diffusive}).}
Define:
$$
\begin{array}{*{20}l}
X := \left[x_1^T,\ldots,x_N^T\right]^T, & S := 1_{N}\otimes s, &
E := X-S,\\
\end{array}
$$
where $1_N$ denotes the $N$-dimensional vector consisting of all 1s. (Notice that such a vector spans $\sM_s$.)
The network dynamics can then be written as:
$$
\dot X = (I_N\otimes A(\sigma(t)))X -(L\otimes \Gamma)X,
$$
so that the error dynamics is described by 
\begin{equation}\label{eqn:error_dyn_net}
\dot E = (I_N\otimes A(\sigma(t)))E -(L \otimes \Gamma)X
\end{equation}
Notice that network coordination is attained if the dynamics of (\ref{eqn:error_dyn_net}) transversal to $\sM_s$ is contracting. Furthermore, notice that
$$
\begin{array}{*{20}l}
(L\otimes \Gamma)X = (L\otimes \Gamma)(E + S) = 
(L\otimes \Gamma)E + (L\otimes \Gamma)S = \\
(L\otimes \Gamma)E + (L\otimes \Gamma)(1_N\otimes s) = 
 (L\otimes \Gamma)E ,
\end{array}
$$
where the last equality follows from Lemma \ref{lem:kronecker} and from the fact that $L \cdot 1_N=0$, since the network is connected by hypothesis. Thus, from (\ref{eqn:error_dyn_net}), we have:
\begin{equation}\label{eqn:error_lin}
\dot E = (I_N\otimes A(\sigma(t))) E - (L\otimes \Gamma) E.
\end{equation}
Since $L$ is symmetric, by means of Lemma \ref{lem:schur} we have that there exist an $N\times N$ orthogonal matrix $Q$ ($Q^TQ=I_N$) such that:
$$
\Lambda=Q^{T}LQ,
$$
where $\Lambda$ is the $N\times N$ diagonal matrix having  the Laplacian eigenvalues as its diagonal elements.

Now, considering the following coordinate transformation: 
$$
Z = \left(Q \otimes I_n\right)^{-1} E,
$$
equation (\ref{eqn:error_lin}) can be recast as
$$
\dot Z= (Q\otimes I_n)^{-1} \left[(I_N\otimes A(\sigma(t))) - (L\otimes \Gamma)\right] (Q\otimes I_n) Z.
$$

Then, using Lemma \ref{lem:kronecker}, we have:
$$
\begin{array}{*{20}l}
\left(Q\otimes I_n\right)^{-1} \left(I_N\otimes A(\sigma(t)) \right) \left(Q\otimes I_n\right) = \\ \left(Q^{-1}\otimes I_n\right) \left(I_N\otimes A(\sigma(t)) \right) \left(Q\otimes I_n\right) = \\
 \left(Q^{-1}\otimes A(\sigma(t))\right)\left(Q\otimes I_n\right) = \\ \left(I_N\otimes A(\sigma(t)) \right).
\end{array}
$$
Analogously:
$$
\begin{array}{*{20}l}
\left(Q\otimes I_n\right)^{-1} \left(L \otimes \Gamma\right)\left(Q\otimes I_n\right)
=\\ \left(Q^{-1}\otimes I_n\right) \left(L \otimes \Gamma\right)\left(Q\otimes I_n\right) = \\
 \left(Q^{-1}L\otimes \Gamma\right)\left(Q\otimes I_n\right) =\\ Q^{-1}LQ\otimes \Gamma = \\ \Lambda\otimes \Gamma.
\end{array}
$$

That is, the network dynamics can be written as:
\begin{equation}\label{eqn:net_lin_transf}
\dot Z = \left[I_N \otimes A(\sigma(t)) - \Lambda \otimes \Gamma\right] Z,
\end{equation}
or equivalently:
$$
\dot z_i =\left[A(\sigma(t)) - \lambda_i\Gamma\right] z_i, \quad i = 1,\ldots, N, \quad z_i \in \R^n.
$$
Now, recall that  the eigenvector associated to the smallest eigenvalue of the Laplacian matrix, i.e. $\lambda_1 = 0$, is $1_N$ and spans $\sM_s$. Therefore, the dynamics along $\sM_s$ is given by
$$
\dot z_1 =\left[A(\sigma(t))\right] z_1,
$$
i.e. it is a solution of the uncoupled nodes' dynamics.
 The dynamics transversal to the invariant subspace is given by:
$$
\dot z_i =\left[A(\sigma(t)) - \lambda_i \Gamma \right] z_i, \quad i = 2,\ldots, N.
$$
Obviously $\left[A(\sigma(t)) - \lambda_i \Gamma \right]$ is bounded and measurable. Thus, by virtue of Corollary \ref{cor:PWL}, all node trajectories globally exponentially converge a.e. towards $\sM_s$, if all of the above dynamics are contracting. 
Now, it is straightforward to check that such a condition is fulfilled if
$$
\dot z_2 =\left[A(\sigma(t)) - \lambda_2 \Gamma \right] z_2
$$ 
is contracting. As this is true from the hypotheses, the result is then proved.

\subsection{A numerical example}
As a representative example, in this Section we use the results presented above to synchronize a network of the form (\ref{eqn:network_dyn}), where the dynamics of each uncoupled node is given by:
\begin{equation}\label{eqn:nodes_dynamics}
\dot x_i:= \left[\begin{array}{*{20}c}
x_{1i}\\
x_{2i}\\
\end{array}\right] = 
\left[\begin{array}{*{20}c}
0 & \abs{\sin(t)}\\
-1 & 0 \\
\end{array}\right] x_i,
\end{equation}
The matrix $\Gamma$  is chosen as:
$$
\Gamma = k \left[\begin{array}{*{20}c}
1 & 0 \\
0 & 1\\
\end{array}\right],
$$
with $k$ being the coupling gain that will be determined using Theorem \ref{thm:linear_diffusive}. The network considered here consists of an all to all topology of three nodes. That is,
$$
L := \left[\begin{array}{*{20}c}
2 & -1 & -1\\
-1 & 2 & -1\\
-1 & -1 & 2\\
\end{array}\right],
$$
and $\lambda_2 = 3$. Thus, from Theorem \ref{thm:linear_diffusive} it follows that the network synchronizes if there exist some matrix measure and two positive scalars $c_1, c_2$ such that:
$$
\begin{array}{*{20}l}
\mu\left(\left[\begin{array}{*{20}c}
0 & \sin(t)\\
-1 & 0 \\
\end{array}\right] - 3 k  \left[\begin{array}{*{20}c}
1 & 0 \\
0 & 1\\
\end{array}\right]\right) \le -c_1,  & \mbox{if}\ \sin(t) \ge 0; \\
\mu\left(\left[\begin{array}{*{20}c}
0 & -\sin(t)\\
-1 & 0 \\
\end{array}\right] - 3 k  \left[\begin{array}{*{20}c}
1 & 0 \\
0 & 1\\
\end{array}\right]\right) \le -c_2,  & \mbox{if}\ \sin(t) < 0.
\end{array}
$$

That is, synchronization is attained if

$$
\begin{array}{*{20}l}
\mu\left(\left[\begin{array}{*{20}c}
-3k & \sin(t)\\
-1 & -3k \\
\end{array}\right] \right) \le -c_1, & \mbox{if}\ \sin(t) \ge 0; \\
\mu\left(\left[\begin{array}{*{20}c}
-3k & -\sin(t)\\
-1 & -3k \\
\end{array}\right] \right) \le -c_2, & \mbox{if}\ \sin(t) <0.
\end{array}
$$
Now, using the matrix measure induced by the vector-$1$ norm (column sums), it is straightforward to see that the above conditions are fulfilled if the coupling gain is selected as
$$
k > \frac{1}{3}.
$$

\begin{figure}%
\includegraphics[width=7cm]{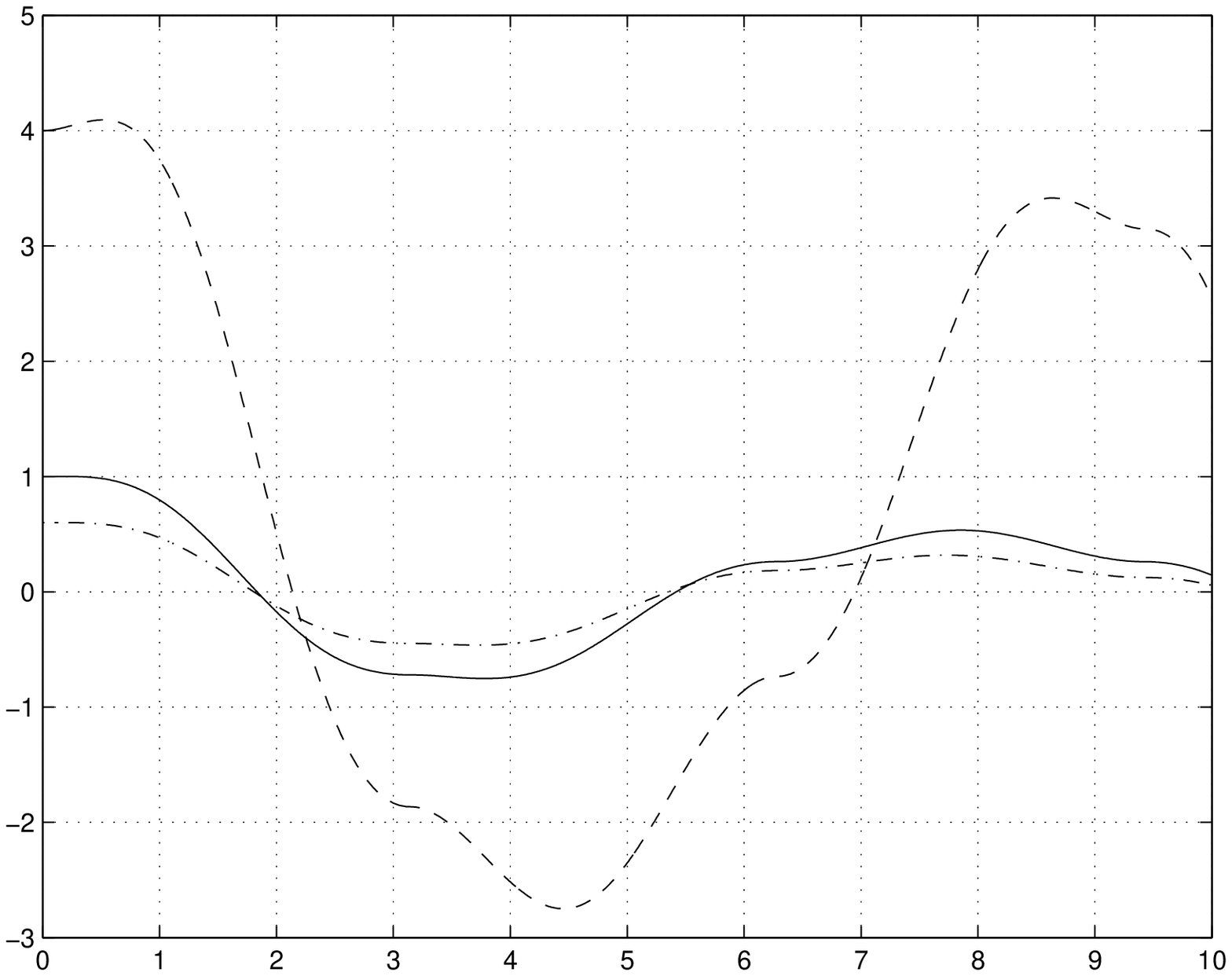}\quad \includegraphics[width=7cm]{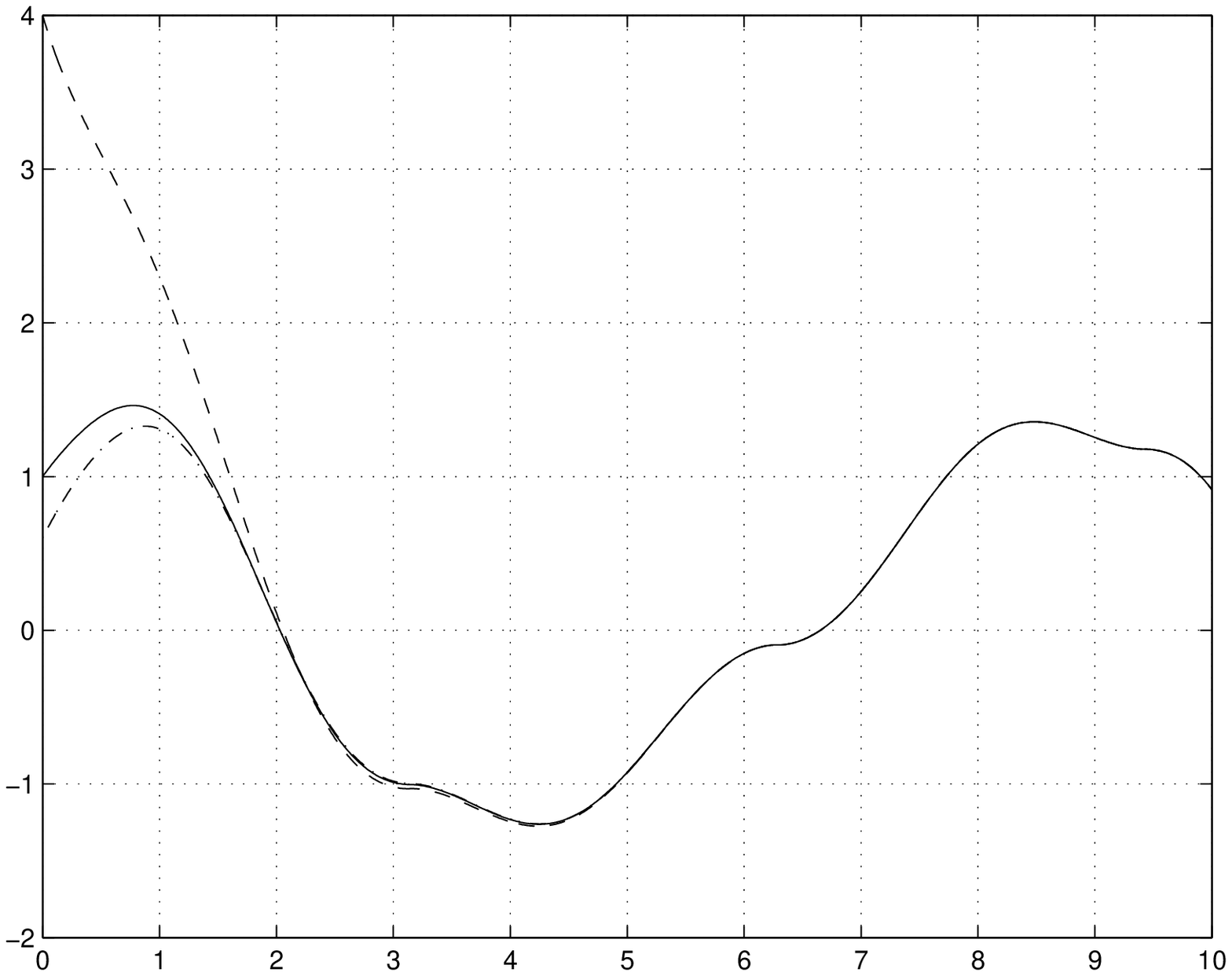}
\caption{State evolution of the network of three switched linear nodes when $k=0$ (left panel) and $k=0.4$ (right panel).}%
\label{fig:threewell-topology}%
\end{figure}

As shown in Fig. \ref{fig:threewell-topology}, the theoretical predictions are confirmed by the numerical simulations.

%

\section{Conclusions}
An extension of contraction theory was presented to a generic class of switched systems: those satisfying conditions for the existence and uniqueness of a Caratheodory solution. In particular, it was proven that infinitesimal contraction of each mode of a switched system of interest gives a sufficient condition for global exponential convergence of trajectories towards each other.   
This result was then used on a set of representative applications. It was shown that, by using contraction, it is possible to immediately derive sufficient conditions for global stability of switched linear systems. Also, contraction was used  to obtain sufficient conditions for the convergence of all nodes in a network of coupled switched linear systems towards the same synchronous evolution.

We wish to emphasize that the results presented in this paper can be immediately applied to generalize to switched Caratheodory systems all of the analysis and design results based on contraction analysis available for smooth systems in the literature. Examples of applications include nonlinear observer design, network protocols design for network coordination, analysis/control of asynchronous systems and biochemical systems.

Finally, note that the results presented in this paper are the first essential stage needed to develop a systematic approach to extend contraction analysis to generic classes of switched systems. The next step is that of addressing the challenging problem of studying convergence in Filippov systems where sliding motion is possible. This is currently work in progress and will be presented elsewhere.


\end{document}